\documentclass[12pt]{amsart}

\newcommand{\C}{\mathbb{C}}
\begin{document}

\title[Bergman kernel: explicit formulas and zeroes]{The Bergman
  kernel function: explicit formulas and zeroes}

\author{Harold P. Boas}
\address{Department of Mathematics\\
Texas A\&M University\\
College Station, Texas, 77843--3368}
\email{boas@math.tamu.edu}

\author{Siqi Fu}
\email{sfu@math.tamu.edu}

\author{Emil J. Straube}
\email{straube@math.tamu.edu}

\thanks{Research supported in part by NSF grant number DMS
  9500916.}

\begin{abstract}
  We show how to compute the Bergman kernel functions of some
  special domains in a simple way.  As an application of the
  explicit formulas, we show that the Bergman kernel functions of
  some convex domains, for instance the domain in~\(\C^3\)
  defined by the inequality \( |z_1|+|z_2|+|z_3|<1\), have
  zeroes.
\end{abstract}

\maketitle

\section{Introduction}
This note has two themes. The first is that one can obtain useful
formulas for the Bergman kernel function in special cases with
minimal computational effort.  The second is that one can find
from the explicit formulas examples of simple convex domains
whose Bergman kernel functions have zeroes.

Recently there has been renewed interest \cite{DAngeloJ1994,
  FrancsicsHanges1996, FrancsicsHanges1997,
  OeljeklausPflugYoussfi} in explicit formulas for the Bergman
kernel function, especially in connection with the asymptotic
behavior of the kernel near weakly pseudoconvex boundary points.
We will show how, even without writing down an infinite series,
one can easily obtain an explicit formula for the Bergman kernel
function of, for example, the domain defined by the inequality \(
\|z_1\|^{2/p_1} + \dots + \|z_n\|^{2/p_n}<1\), where each~\(z_j\)
is a vector in some \(\C^{m_j}\), and each~\(p_j\) is a positive
integer.

Finding a simple class of domains whose Bergman kernel functions
are zero-free has been an open problem ever since Lu Qi-Keng
\cite{LuQiKeng1966} raised the question in connection with the
existence of Bergman representative coordinates.  The first
author recently discovered \cite{BoasH1996} that the Bergman
kernel of a generic strongly pseudoconvex domain does have
zeroes, but the possibility was left open that all convex domains
might have zero-free Bergman kernel functions. We will show that
the Bergman kernel function of the domain in~\(\C^3\) defined by
the inequality \(|z_1|+|z_2|+|z_3|<1\) does have zeroes.
Consequently, when \(n\ge 3\), there exists a smooth, bounded,
strongly convex Reinhardt domain in \(\C^n\), even with real
analytic boundary, whose Bergman kernel function has zeroes.

\section{Explicit formulas for the Bergman kernel}
\label{sec:explicit}
The Bergman kernel function \(K(z,w)\) of a domain~\(\Omega\) in
\(\C^n\) is the reproducing kernel for the space of
square-integrable holomorphic functions, that is, \(
f(z)=\int_\Omega f(w)K(z,w)\,dV_w\) when \(f\)~is a
square-integrable holomorphic function in~\(\Omega\).  We use
three basic principles to obtain new Bergman kernel functions
from old ones. The principles of deflation and inflation relate
kernel functions of domains in different dimensions. On the other
hand, the well known principle of folding relates kernel
functions of domains in the same dimension.

\subsection{Deflation}
Let \(\Omega\) be a bounded domain in~\(\C^n\) defined by an
inequality of the form \( \varphi(z)<1\), where \(\varphi\)~is a
continuous, nonnegative function in a neighborhood of the closure
of~\(\Omega\). Let \(K_1\) denote the Bergman kernel function of
the bounded domain~\(\Omega_1\) in \(\C^{n+2}\) defined by the
inequality \(\varphi(z) + |\zeta_{1}|^{2/p}+
|\zeta_{2}|^{2/q}<1\), where \(p\) and~\(q\) are positive real
numbers. Similarly, let \(K_2\) denote the Bergman kernel
function of the bounded domain~\(\Omega_2\) in \(\C^{n+1}\)
defined by the inequality \(\varphi(z)+ |\zeta|^{2/(p+q)}<1\).
Then
\begin{equation}
  \label{deflation}
  \pi K_2(z,0, w,0) = \frac{\pi^2
    \Gamma(p+1)\Gamma(q+1)}{\Gamma(p+q+1)} K_1(z, 0,0, w, 0, 0).
\end{equation}

For example, we will show in~\S\ref{sec:two dim} that the Bergman
kernel function \(K(z_1,z_2,w_1,w_2)\) of the non-convex domain
\( \{ (z_1,z_2): |z_1|+|z_2|^{2/4}<1\}\) in~\(\C^2\) has zeroes
when \(z_2=w_2=0\).  Taking \(n=1\) and \(p=q=2\)
in~\eqref{deflation} shows that the Bergman kernel function of
the convex domain \( \{ (z_1,z_2,z_3): |z_1|+|z_2|+|z_3|<1\}\)
in~\(\C^3\) also has zeroes.

The deflation identity~\eqref{deflation} holds because both sides
represent the (unique) holomorphic reproducing kernel function
for the holomorphic functions on~\(\Omega\) that are
square-integrable with respect to the weight function \(
(1-\varphi)^{p+q}\).  (Consult \cite{LigockaE1989} for a
systematic study of weighted Bergman projections.)  To see this,
observe that a holomorphic function \(h\) on~\(\Omega\) can be
regarded as a holomorphic function on both \(\Omega_1\)
and~\(\Omega_2\) by extending it to be independent of the extra
variables. The reproducing property of the Bergman kernel
function on~\(\Omega_2\) implies that
\begin{equation}
\label{reproducing}
h(z)= \int_{\Omega_2} h(w) K_2(z,0, w, \eta) dV_{w,\eta}.
\end{equation}
In the domain~\(\Omega_2\), the fiber over a point~\(w\)
in~\(\Omega\) is a one-dimensional disc of radius \(
(1-\varphi(w) )^{(p+q)/2}\).  Integrating with respect
to~\(\eta\) and using the mean-value property of harmonic
functions shows that~\eqref{reproducing} reduces to
\begin{equation*}
  h(z)= \int_\Omega h(w) (1-\varphi(w) )^{(p+q)} \pi
  K_2(z,0,w,0)\,dV_{w}.
\end{equation*}
The subaveraging property of \( |K_2(z,0,w,0)|^2\) implies that
\(K_2(z,0,w,0)\) is square-integrable with respect to the weight
\((1-\varphi)^{p+q}\), so the left-hand side of~\eqref{deflation}
is the reproducing kernel for the holomorphic functions
on~\(\Omega\) that are square-integrable with respect to the
weight \((1-\varphi)^{p+q}\).

The verification of the reproducing property of the right-hand
side of~\eqref{deflation} differs only in that the fiber over a
point in~\(\Omega\) is two-dimensional, and one needs the volume
of the region in~\(\C^2\) defined by the inequality \(
|\zeta_1|^{2/p}+|\zeta_2|^{2/q}<R\).  A routine calculation in
polar coordinates shows that this volume equals \( \pi^2
\Gamma(p+1)\Gamma(q+1) R^{p+q}/ \Gamma(p+q+1)\), which completes
the proof of the deflation identity~\eqref{deflation}.

\subsection{Inflation}
\label{sec:inflation}
Let \(\Omega\) be a bounded complete Hartogs domain
in~\(\C^{n+1}\) defined by an inequality of the form \( |\zeta|^2
<\varphi(z)\), where \(\zeta\in\C\), \(z\in\C^n\), and
\(\varphi\) is a bounded, positive, continuous function on the
interior of some bounded domain in~\(\C^n\). Due to the circular
symmetry in the one-dimensional variable, there is a function
\(L(z,w,t)\) such that the Bergman kernel function
\(K(z,\zeta,w,\eta)\) for~\(\Omega\) can be written in the form
\(L(z,w,\zeta\bar\eta)\).  We can inflate \(\Omega\) into a
domain~\(\widetilde\Omega\) in \(C^{n+m}\) defined by the
inequality \( \|Z\|^2<\varphi(z)\), where \(Z\)~is a vector
variable in~\(\C^m\) with Euclidean norm~\(\|Z\|\).  We will show
that the Bergman kernel function \(\widetilde K(z,Z,w,W)\)
of~\(\widetilde \Omega\) is given by the relation
\begin{equation}
\label{inflation}
\widetilde K(z,Z,w,W) = \frac{1}{\pi^{m-1}}
\left.\frac{\partial^{m-1}}{\partial t^{m-1}}
  L(z,w,t)\right|_{t=\langle Z,W\rangle},
\end{equation}
where \(\langle Z,W\rangle= Z_1\overline W_1+ \dots+ Z_m\overline
W_m\).

The motivation for this formula is the relation between the
Bergman kernel function \(\pi^{-1}(1-\zeta\bar\eta)^{-2}\) for
the unit disc in~\(\C\) and the Bergman kernel function for the
unit ball in~\(\C^m\), which is
\begin{equation}
\label{eq:ball}
  \frac{1}{\pi^{m}} \left.\frac{\partial^{m-1}}{\partial
        t^{m-1}} \frac{1}{(1-t)^2}\right|_{t=
      \langle Z,W \rangle} = \frac{m!}{\pi^m}\cdot
    \frac{1}{(1-\langle Z,W \rangle)^{m+1}}.
\end{equation}
For an application, recall that when \(p\)~is a positive real
number, the Bergman kernel function of the domain
\(|z|^2+|\zeta|^{2/p}<1\) in~\(\C^2\) is
\begin{equation}
\label{eq:Bergman}
  \frac{1}{p\pi^2} \left. \frac{\partial^2}{\partial t^2}
  \left(\frac{1}{(1-t)^p-\zeta\bar\eta}\right) \right|_{t=z\bar w}.
\end{equation}
(Here is the only place in this note where---implicitly---an
infinite series appears.  Bergman \cite[p.~82]{BergmanS1936}
computed this kernel by summing a series; although he states that
\(p\)~is integral, his computation is valid for arbitrary
positive~\(p\).)  Inflating both variables shows that the Bergman
kernel function of the domain \(\|z\|^2+\|Z\|^{2/p}<1\)
in~\(\C^{n+m}\) is
\begin{equation}
\label{eq:pflate}
\frac{1}{p\pi^{n+m}} \left.\frac{\partial^{n+m}}{\partial
    t^{n+1}\partial u^{m-1}} \left(\frac{1}{(1-t)^p-u}\right)
\right|_{\substack {t=\langle z,w\rangle \\ u=\langle
    Z,W\rangle}}.
\end{equation}
This kernel function was first computed (in a somewhat different
form) by Chalmers \cite[Theorem~1.1]{ChalmersBL1969}.
Subsequently, D'Angelo used a different method to sum the
orthonormal series to compute this kernel for \(m=1\) in
\cite{DAngeloJ1978} and for general~\(m\) in \cite{DAngeloJ1994};
his formulas are written as finite sums rather than as
derivatives.

We prove the inflation formula~\eqref{inflation} by verifying the
reproducing property.  Since \(\widetilde\Omega\)~is invariant
under unitary transformations of~\(\C^m\), it suffices to
verify~\eqref{inflation} when \(Z=(Z_1,0,\dots,0)\), in which
case we can rewrite
\begin{equation}
\label{eq:rewrite}
  \left.\frac{\partial^{m-1}}{\partial t^{m-1}}
    L(z,w,t)\right|_{t=\langle Z,W\rangle} =
  \frac{1}{Z_1^{m-1}} \frac{\partial^{m-1}}{\partial
        \overline W_1^{m-1}} L(z,w,Z_1\overline W_1).
\end{equation}

If \(h\)~is a square-integrable holomorphic function
on~\(\widetilde\Omega\), then we find by using the mean-value
property in the \(m-1\) variables \(W_2\), \dots, \(W_m\) that
\begin{multline*}
  \frac{1}{\pi^{m-1}} \int_{\widetilde\Omega} h(w,W) 
\left.\frac{\partial^{m-1}}{\partial t^{m-1}}
  L(z,w,t)\right|_{t=Z_1\overline W_1} \,dV_{w,W} 
= \frac{1}{(m-1)!\,Z_1^{m-1}} \times \\
\times \int_{\Omega} h(w,W_1) 
(\varphi(w) -|W_1|^2)^{m-1}
\frac{\partial^{m-1}}{\partial\overline W_1^{m-1}}
L(z,w,Z_1\overline W_1)\,dV_{w,W_1}.
\end{multline*}
Observe that \( f(z) \mapsto (1-|z|)f'(z)\) is continuous
in~\(L^2\) on holomorphic functions in the unit disc (this
follows by integrating in polar coordinates), so the above
calculation implies that the kernel~\eqref{eq:rewrite} is
square-integrable in~\(\widetilde\Omega\).  Approximating
\(h(w,W_1)\) and \( L(z,w,Z_1\overline W_1)\) by polynomials in
\(W_1\) and~\(\overline W_1\) respectively shows that no boundary
terms arise when we integrate by parts in the \(W_1\)~variable
\((m-1)\) times, so the right-hand side of the preceding formula
reduces to
\begin{equation*}
  \frac{1}{Z_1^{m-1}}\int_{\Omega} h(w,W_1) W_1^{m-1}
  L(z,w,Z_1\overline W_1)\,dV_{w,W_1} = h(z, Z_1),
\end{equation*}
where the final equality follows from the reproducing property of
the Bergman kernel function on~\(\Omega\).  Thus the kernel given
by~\eqref{inflation} does reproduce holomorphic functions
on~\(\widetilde\Omega\).

\subsection{Folding}
The third method is a well known ``useful tool to establish
explicit formulas for the Bergman kernel'' (Jarnicki and Pflug
\cite[p.~176]{JarnickiPflug1993}).  As in
paragraph~\ref{sec:inflation}, let \(\Omega\) be a bounded
complete Hartogs domain in \(\C^{n+1}\) defined by an inequality
of the form \( |\zeta|^2 <\varphi(z)\).  When \(p\)~is a positive
integer, \( (z,\zeta)\mapsto (z,\zeta^p)\) is a proper
holomorphic mapping of~\(\Omega\) onto the domain~\(\Omega_p\) in
\(\C^{n+1}\) defined by the inequality \(
|\zeta|^{2/p}<\varphi(z)\).  Bell's transformation rule
\cite{BellS1982} for the Bergman kernel function under proper
mappings implies that when \(\zeta\bar\eta\ne 0\), the Bergman
kernel~\(K\) for~\(\Omega\) is related to the Bergman
kernel~\(K_p\) for~\(\Omega_p\) via
\begin{equation*}
  p^2 (\zeta\bar\eta)^{p-1} K_p(z,\zeta^p,w,\eta^p) =
  \sum_{j=1}^p \bar\omega^j K(z,\zeta,w,\omega^j \eta),
\end{equation*}
where \(\omega\)~is a primitive \(p\)th root of unity.

For example, we can repeatedly fold formula~\eqref{eq:pflate}
with \(m=1\) to see that when \(p_1\), \dots, \(p_n\) are
positive integers, and \(p\)~is any positive real number, the
domain in \(\C^{n+1}\) defined by the inequality \( |z_1|^{2/p_1}
+ \dots + |z_n|^{2/p_n} +|z_{n+1}|^{2/p}<1\) has Bergman
kernel~\(K\) given by the finite sum
\begin{multline}
\label{eq:general}
K(z_1^{p_1},\dots, z_n^{p_n}, z_{n+1}, w_1^{p_1}, \dots,
w_n^{p_n}, w_{n+1}) = \frac{1}{p\pi^{n+1}} \prod_{k=1}^n
\frac{1}{ p_k^2 (z_k\bar w_k)^{p_k-1}} \times \\ \times
\sum_{j_1=1}^{p_1} \dots\sum_{j_n=1}^{p_n}
\frac{\partial^{n+1}}{\partial t^{n+1}} \left.  \left(
    \frac{\bar\omega_1^{j_1}\dots\bar\omega_n^{j_n}} {(1-t)^p
      -z_{n+1} \bar w_{n+1}} \right) \right|_{t=\sum_{\ell=1}^n
  z_\ell \bar w_\ell \bar \omega_\ell^{j_\ell}},
\end{multline}
where each \(\omega_j\) is a primitive \(p_j\)th root of unity.
Francsics and Hanges \cite[Theorem~3]{FrancsicsHanges1996}
computed this kernel by summing an infinite series. When \(p\)~is
also an integer, this kernel can be obtained by folding the
simpler formula~\eqref{eq:ball} for the ball (compare
\cite[Exercise~6.5]{JarnickiPflug1993}); this special case was
first computed by Zinov$'$ev \cite{ZinovevBS1974} via series
methods.

It should be evident how to inflate~\eqref{eq:general}
via~\eqref{inflation} to obtain the Bergman kernel function for
the domain defined by the inequality \( \|z_1\|^{2/p_1}+\dots +
\|z_n\|^{2/p_n} +\|z_{n+1}\|^{2/p}<1\), where each \(z_j\) is a
vector in some~\(\C^{m_j}\).  Repeated foldings and inflations
easily give the Bergman kernel functions of still more general
domains with defining functions like \( (\|z_1\|^{2/p_1} +
\|z_2\|^{2/p_2})^{2/p} + \|z_3\|^{2/p_3}<1\), the only difficulty
in writing the formulas being the choice of a suitable notation
(compare \cite[\S6.2]{EgorychevGP1984}).

\section{Zeroes of the Bergman kernel}
We now deduce from the explicit formulas some consequences about
when the Bergman kernel function has zeroes. In particular, the
kernel function has zeroes in some simple convex domains, a fact
previously unknown.

\subsection{Two-dimensional domains}
\label{sec:two dim}
To find the Bergman kernel function~\(K_p\) for the
domain~\(\Omega_p\) in~\(\C^2\) defined by the inequality \(
|z_1| + |z_2|^{2/p}<1\), we fold Bergman's
formula~\eqref{eq:Bergman}, or, equivalently, take \(n=1\) and
\(p_1=2\) in~\eqref{eq:general}, leaving~\(p\) as an arbitrary
positive real number.  Setting \(x:=z_1\bar w_1\) and
\(y:=z_2\bar w_2\), we have
\begin{equation}
\label{eq:bergman}
  K_p(z_1^2,z_2,w_1^2,w_2)= \frac{1}{4p\pi^2 x} \cdot
  \frac{\partial^2}{\partial x^2} \left( \frac{1}
    {(1-x)^p-y}-\frac{1} {(1+x)^p-y} \right)
\end{equation}
when \(x\ne0\).  We will show that the Bergman kernel
function~\(K_p\) of~\(\Omega_p\) has zeroes when \(p>2\), while
the Bergman kernel function~\(K_2\) of~\(\Omega_2\) has no zeroes
on the interior of the domain.

Setting \(y=0\) in~\eqref{eq:bergman}, we see that
\(K_p(z_1^2,0,w_1^2,0)=0\) if and only if \(x\ne0\) and
\((1+x)^{p+2}=(1-x)^{p+2}\).  When \(y=0\), the point~\(x\) runs
over the unit disc. Since \( t\mapsto (1+t)/(1-t)\) maps the unit
disc bijectively to the right half-plane, it follows that the
equation \([(1+x)/(1-x)]^{p+2}=1\) has a non-zero solution in the
unit disc if and only if \(p>2\). Hence for \(p>2\), the Bergman
kernel function~\(K_p\) has zeroes when \(z_2=0\).

Taking the limit in~\eqref{eq:bergman} as \(x\to0\), we find
after a routine calculation that
\begin{equation*}
  K_p(0,z_2,0,w_2)= \frac{1}{2\pi^2} \cdot \frac{y^2(p^2-3p+2) +
    4y(p^2-1) +(p^2+3p+2)}{(1-y)^4}.
\end{equation*}
This expression is real-valued for real values of~\(y\), positive
for \(y=0\), and equal to \(-(p^2-4)/\pi^2\) at \(y=-1\), so if
\(p>2\), then there is a zero for \(y\)~in the unit disc. Hence
for \(p>2\), the Bergman kernel~\(K_p\) has zeroes also when
\(z_1=0\).

Setting \(p=2\) in~\eqref{eq:bergman}, we find after some routine
algebra that
\begin{equation*}
  K_2(z_1,z_2,w_1,w_2) = \frac{2}{\pi^2}\cdot
  \frac{3(1-x-y)(1-(x-y)^2) +8xy}{ ((1-x-y)^2-4xy)^3},
\end{equation*}
where \(x=z_1\bar w_1\) and \(y=z_2\bar w_2\) as before (compare
\cite[Example 6.1.9]{JarnickiPflug1993}). We aim to show that
this expression has no zeroes when \(x\) and~\(y\) are complex
numbers such that \(\sqrt{|x|}+\sqrt{|y|}<1\). This constraint
implies that \(4|x||y|<(1-|x|-|y|)^2\), so the numerator of the
expression exceeds twice
\begin{multline*}
  3(1-|x|-|y|)(1-|x-y|^2)-2(1-|x|-|y|)^2\\
\ge 3(1-|x|-|y|)^2(1+|x-y|)-2(1-|x|-|y|)^2\\
> (1-|x|-|y|)^2.
\end{multline*}
Consequently, \(K_2(z_1,z_2,w_1,w_2)\) has no zeroes on the
interior of the domain. However, there are zeroes on the
boundary: for instance, \(K_2(1,0,-1,0)=0\).

Notice that when \(p>2\), the case in which we have shown that
the Bergman kernel~\(K_p\) has zeroes, the domain~\(\Omega_p\) is
not convex. In higher dimensions, the situation changes.

\subsection{Higher-dimensional domains}
Using the deflation identity~\eqref{deflation}, we can identify
many convex domains in higher dimensions whose Bergman kernel
functions have zeroes. For example, the Bergman kernel function
of the convex domain in~\(\C^n\) defined by the inequality
\(|z_1|+|z_2|+\dots+|z_n|<1\) has zeroes if and only if the
dimension \(n\ge 3\). Indeed, we have checked that this Bergman
kernel is zero-free when \(n=2\). If \(n\ge 3\), on the other
hand, then we find by~\eqref{deflation} that the restriction of
the kernel to the subspace where \(z_2=\dots=z_n=0\) equals a
positive constant times \(K_{2n-2}(z_1,0,w_1,0)\), and we saw
in~\S\ref{sec:two dim} that this kernel has zeroes when
\(2n-2>2\).

As another example, we show that the Bergman kernel
function~\(K\) of the convex domain in~\(\C^n\) defined by the
inequality \(|z_1|+|z_2|^2+|z_3|^2+\dots+|z_n|^2<1\) has zeroes
if and only if the dimension \(n\ge 4\).  Using deflation, as in
the preceding paragraph, shows that the kernel has zeroes when
\(n\ge 4\), but does not show that the kernel is zero-free when
\(n\le 3\).  Accordingly, we examine the explicit formula
obtained by folding the kernel for the ball:
\begin{multline*}
  K(z_1^2,z_2,\dots,z_n,w_1^2, w_2, \dots, w_n)\\ =
  \frac{n!}{\pi^n} \cdot \frac{1}{4z_1\bar w_1} \left(
    \frac{1}{(1-\langle z,w\rangle)^{n+1}} - \frac{1}{(1+\langle
      z,w\rangle)^{n+1}} \right).
\end{multline*}
The same argument as in \S\ref{sec:two dim}, that \( t\mapsto
(1+t)/(1-t)\) maps the unit disc bijectively to the right
half-plane, shows that this kernel function has zeroes in the
interior of the domain if and only if \(n+1>4\).

By using the tools indicated in~\S\ref{sec:explicit}, the reader
can easily produce more examples of convex domains whose Bergman
kernel functions have zeroes. We have given two examples to
illustrate the usefulness of the techniques.

The examples just given of convex domains whose Bergman kernel
functions have zeroes are nonsmooth domains.  We can exhaust such
domains by sequences of smooth, strongly convex domains with
real-analytic boundaries. By Ramadanov's theorem
\cite{RamadanovI1967}, the Bergman kernels of the approximating
domains converge uniformly on compact sets to the Bergman kernel
of the limit domain. By Hurwitz's theorem, the Bergman kernel
functions of the approximating strongly convex domains cannot all
be zero-free.  Consequently, when \(n\ge 3\), there exist smooth,
bounded, strongly convex domains in~\(\C^n\), even with
real-analytic boundary, whose Bergman kernel functions have
zeroes.

\section{Remarks}
\begin{enumerate}
\item The domains in~\(\C^n\) whose Bergman kernels we have shown
  to have zeroes are convex only when \(n\ge 3\).  It would be
  interesting to find a bounded convex domain in dimension
  \(n=2\) whose Bergman kernel function has zeroes.

\item It would be interesting to characterize the set of positive
  real numbers \(p_1\), \(p_2\), \dots, \(p_n\) for which the
  Bergman kernel function of the domain in~\(\C^n\) defined by
  the inequality \( |z_1|^{2/p_1} + |z_2|^{2/p_2}+ \dots +
  |z_n|^{2/p_n}<1\) is zero-free.

\end{enumerate}

\subsection*{Note added 23 June 1997}
Shortly after we emailed P.~Pflug about our results, he and E.~H.
Youssfi were able to find zeroes of the Bergman kernels of the
``minimal balls'' studied in \cite{OeljeklausPflugYoussfi}.

\newpage           


\providecommand{\bysame}{\leavevmode\hbox to3em{\hrulefill}\thinspace}

\end{document}